\documentclass[12pt]{article}
\usepackage{latexsym}
\usepackage{amssymb}

\newcommand{\proof}{{\noindent \bf Proof. }}
\newtheorem{thm}{Theorem}

\newtheorem{prop}{Proposition}

\newcommand{\D}{{\cal D}}
\newcommand{\F}{{\cal F}}
\newcommand{\G}{{\cal G}}

\date{}

\begin{document}
\begin{titlepage}
\title{\bf FAMILIES OF LOCALLY SEPARATED HAMILTON PATHS}
{\author{{\bf  J\'anos K\"orner}
\\e--mail: {\tt korner@di.uniroma1.it} 
\and{\bf Angelo Monti}
\\e--mail:  {\tt <monti@di.uniroma1.it>}\quad
\medskip
\\Sapienza University of Rome\\
\\ ITALY}}

\maketitle
\begin{abstract}
We improve by an exponential factor the lower bound of K\"orner and Muzi for the cardinality of the largest family of Hamilton paths in a complete graph of 
$n$ vertices in which the union of any two paths has degree 4. The improvement is through an explicit construction while the previous bound was obtained by a greedy algorithm. We solve a similar problem for permutations up to an exponential factor.  
\bigskip

{\bf Keywords} Hamilton paths, graph-difference, permutations, intersection problems

{\bf AMS Subject classification numbers} 05D99, 05C35, 05C62, 94A24

\end{abstract}
\end{titlepage}

\section{Introduction}

Greedy algorithms play a central role in combinatorial problems related to information theory. The classical example for this is the code distance problem asking for the maximum cardinality of a set of binary sequences of length $n$ any two of which differ in at least $\alpha n$ coordinates  where $\alpha$ is a fixed constant. Gilbert \cite{Gi} and Varshamov \cite{V} showed independently that for $\alpha \in [0, 1/2)$ the maximum cardinality is at least $\exp_2{n(1-h(\alpha))}$ where $h$ is the binary entropy function 
$h(t)=-t \log_2-(1-t) \log_2(1-t).$ This result, called the Gilbert--Varshamov bound, is proved by a greedy algorithm. Considerable effort went into trying to strengthen this bound, 
and although it was recently improved by a linear factor \cite{Var}, \cite{Vu}, it is believed by many to be exponentially tight. 

In an attempt to generalise Shannon's zero--error capacity problem \cite{Sh} from fixed--length sequences of elements of a finite alphabet to permutations and then to 
subgraphs of a fixed complete graph, cf. \cite{KM}, \cite{sev} and \cite{KMM}, one of the main construction techniques is again based on greedy algorithms. It is crucial to understand the true role of greedy algorithms in our context. We will study two different combinatorial problems of information--theoretic flavour and discuss the constructions 
obtained by the greedy algorithm in both cases. 
 For the general framework and especially for Shannon's zero--error capacity problem we refer the reader to 
the monograph \cite{CK}.
\section{Hamilton paths}
The first problem we are to deal with is about Hamilton paths in the fixed complete graph $K_n$ with vertex set $[n]=\{1,2,\dots, n\}.$ We will say that two such paths are 
{\em} crossing if their union has degree 4. Let us denote by $Q(n)$ the maximum cardinality of a family of Hamilton paths from $K_n$ in which any pair of paths is 
crossing. It was proved in \cite{KMM} (in Theorem 3) that
$$ \frac{(n-2)!}{\lfloor n/2 \rfloor ! (1+\sqrt{2})^n}\leq Q(n)\leq \frac{n!}{\lfloor n/2 \rfloor!  2^{\lfloor n/2 \rfloor}}.$$
The lower bound was achieved by a greedy algorithm. The core of this proof was an upper bound on the number of Hamilton paths not crossing a fixed path. 
Here we give an explicit construction resulting in a larger lower bound. The main message of this result is that for the present problem the greedy algorithm is not tight. The improvement is in the second term of the asymptotics of $Q(n)$ and it could not be otherwise for the leading terms in the lower and upper bounds for $Q(n)$ in \cite{KMM} 
coincide. In fact, both of the bounds in the above theorem are of the form 
$2^{cn} (dn)!$
where $c$ and $d$ are constants. The leading term is the factorial and the constant $d$ in it is the same in both of the bounds. The difference is in the constant in the exponential factor. Our new result increases the constant in the exponential factor in the lower bound in Theorem 3 in \cite{KMM}.

More precisely, we have

\begin{thm} \label{thm:sina} 
$$\frac{(\lfloor n/2 \rfloor -1)!}{2^{\lfloor n/4 \rfloor }}\leq Q(n).$$

\end{thm}

Before proving our result we observe that
$$ \frac{(n-2)!}{\lfloor n/2 \rfloor ! (1+\sqrt{2})^n}  < \Big(\frac{2}{1+\sqrt{2}}\Big)^n\cdot \lceil n/2 \rceil ! \leq\frac{(\lfloor n/2 \rfloor-1) !}{2^{\lfloor n/4 \rfloor }}$$
where the first inequality follows since 
$$\frac{(n-2)!}{\lfloor n/2 \rfloor !}\leq \frac{n!}{\lfloor n/2 \rfloor !}={n \choose {\lfloor n/2 \rfloor }}( \lceil n/2 \rceil)!<2^n( \lceil n/2 \rceil)!$$
while the last inequality is obvious as $\frac{2}{1+\sqrt{2}} <  2^{-1/4}.$ 
Thus the new lower bound of this theorem improves on the one in \cite{KMM} quoted above.

\medskip

\proof

We consider the complete bipartite graph $K_{\lfloor n/2 \rfloor , \lceil n/2 \rceil}.$ For  the ease of notation let us suppose that this graph is defined by the bipartition 
of $[n]$ into the sets $A$ and $B$ where $A=[\lfloor n/2 \rfloor]$ is the set of the first  $\lfloor n/2 \rfloor$ natural numbers while $B=[n]-A.$ Let us fix a permutation of the elements of $B.$ Then any permutation of $A$ defines a Hamilton path in $K_{\lfloor n/2 \rfloor , \lceil n/2 \rceil}$ where the first element of the path is the first element of the permutation of $B$. The path alternates the elements of $A$ and $B$ in such a way that the elements of $B$ are in the odd positions along the path and in the order defined by the 
permutation of $B$ we have fixed. Likewise, the relative order of the elements of $A$ following the path is the one given by the permutation of $A.$ Every 
Hamilton path in this family can thus be identified with some permutation of the set $A.$ We represent each permutation as a linear order of $A.$

Notice now that if there is an element $a \in A$ such that in two different linear orders its position differs by at least 2 while none of these positions is the last one, then the 
vertex $a$ has disjoint pairs of adjacent vertices in the two Hamilton paths. In other words, the vertex $a$ has degree 4 in the union of these paths. We will say that such a pair of permutations is 2-different. We will now construct a sufficiently large family of pairwise 2--different permutations. To this end, we first observe that if two permutations are 
2--different, then in their respective inverses there is some natural number whose images differ by at least two. K\"orner, Simonyi and Sinaimeri have shown 
(cf. Theorem 1 in \cite{KSS}) that for any natural number $m$ the set $[m]$ has a family of exactly $\frac{m!}{2^{\lfloor m/2 \rfloor}}$ permutations such that for any two of them there is at least one element of $[m]$ that has images differing by at least two. Hence the inverses of these permutations form a family satisfying our original condition, except for the fact that the $a$ whose positions in a pair of permutations differ by at most two positions might be in the last position in one of them.  To exclude this event, we choose a subfamily of our permutations in which the last element in all the permutations is the same.  Choosing $m=|A|$ we complete the proof.

\hfill$\Box$

The construction in the last proof seems sufficiently natural for us to believe that it might be asymptotically optimal at least for Hamilton paths in a balanced complete 
bipartite graph. In other words, let $B(n)$ denote the largest cardinality of a family of Hamilton paths in $K_{\lfloor n/2 \rfloor , \lceil n/2 \rceil}$ any two members of 
which are crossing. If the answer is positive, one might go even further and ask whether the lower bound in Theorem \ref{thm:sina} is asymptotically tight.

\section{Intersecting vs. separated families}

In \cite{KMM} where the previous problem was introduced, it was presented alongside with a fairly general framework for extremal combinatorics with an information-theoretic flavour. This framework is in terms of  two (mostly but not necessarily disjoint) families of graphs on the same vertex set $[n]=\{1,2,\dots, n\}$. If we denote these families by 
$\F$ and $\D,$ respectively, then
$M(\F,\D, n)$ stands for the largest cardinality of 
a subfamily $\G\subseteq \F$ with the property that the union of any two of its different members belongs to $\D$. Here the union of two graphs on the same vertex set is the graph whose edge set is the union of those of the two graphs. This framework is the fruit of an attempt to describe a consistent part of extremal combinatorics where information theoretic methods seem to 
be relevant. In \cite{KMM} it was emphasised that when the families $\F$ and $\D$ are disjoint, the question about $M(\F,\D, n)$ is fairly information--theoretic in nature and there never is a natural candidate for an optimal construction; much as it happens for the graph capacity problem of Shannon's. As a matter of fact, in each of these cases 
one is interested in the largest cardinality of a family of mathematical objects (graphs or strings) any pair of which are different in some specific way. Intuitively, the objects in these families are not only different, but also distinguishable in some well--defined sense inherent to the problem. In other words, the optimal constructions are, in a vague and broad sense, similar to the channel codes in information theory.

If, however, $\F=\D,$ the solution of our problems seems to be of a completely different nature, and, in particular, any analogy with codes disappears. Rather, the optimal solutions 
often are so--called kernel structures, and even when this is not the case, kernel structures give good constructions for the problem. Kernel structures are families of objects having 
in common a fixed projection, called the kernel. This means that in a vague and broad sense, these objects are similar.
We are trying to build up a dichotomy distinguishing those problems for which kernel structures can be defined (and often are even optimal) and those where the optimal structures are like "codes" in information theory.  It would be highly interesting to get a better understanding of this dichotomy. 

At the surface, the dichotomy for the problems of determining $M(\F,\D, n)$ is based on whether or not $\F=\D.$
However, this view is too simplistic. As a matter of fact, 
our next problem will illustrate the lack of an easy way to distinguish between the two cases. More precisely, we will give an example where the families $\F$ and $\D$ are disjoint, and the  solution of the problem still has a kernel structure. To explain our example, 
let $\F_{\rm cy}$ be the family of Hamilton cycles in the complete graph $K_n.$ However, let $\D_3$ be the family of graphs having at least one vertex of 
degree 3. The relation of having a vertex of degree 3 in the union of two graphs is clearly irreflexive when restricted to cycles. In fact, the problem of determining 
$M(\F_{\rm cy}, \D_3,n)$ is, on the surface, no different from the one which is the subject of Theorem \ref{thm:sina}. However, this impression is wrong and the present problem is completely different from the previous one. 
To understand why, it suffices to realise that the union of two Hamilton cycles in $K_n$ has a vertex of degree 3 if and only if they share a common edge. 
This makes it possible to invoke kernel structures. Moreover, not surprisingly, a quasi--optimal solution for this problem has a kernel structure.  In fact, we have 

\begin{prop}\label{prop:cy}

$$M(\F_{\rm cy}, \D_3,n)\geq(n-2)!$$
and this bound is tight for odd $n.$ For even $n$ we have
$$M(\F_{\rm cy}, \D_3,n)\leq (n-1)[(n-3)!]$$
\end{prop}

\proof

To prove the lower bound, consider the family of all Hamilton cycles in $K_n$ containing a fixed edge. To establish the upper bound, let first $n$ be odd. By a classical construction of Walecki (1890), cited in \cite{Ak}, for every odd $n$, the family of the Hamilton cycles of $K_n$ can be partitioned into $(n-2)!$ classes where the cycles belonging to the same class are disjoint. (Each class is a partition of the edges of $K_n$ into $\frac{n-1}{2}$ Hamilton cycles). This gives our upper bound for odd $n$. For even $n$ we can only have families of $\frac{n-2}{2}$ pairwise disjoint Hamilton cycles which gives the weaker bound for this case.

\hfill$\Box$

\section{Permutations}

A permutation of $[n]$ is nothing but a (consecutively oriented) Hamilton path in the symmetrically complete directed graph $K_n.$ (With a slight abuse of notation we denote by $K_n$ also the directed graph in which every pair of distinct vertices is connected by two edges going in opposite direction.) In the previous part of the paper we were dealing with problems about Hamilton paths in (the undirected complete graph) $K_n$ such that for any pair of them, for at least one vertex in $K_n$ the union of its 1-neighbourhoods in the two paths satisfied some disjointness condition. We now turn to an analogous problem about directed paths. 

In a directed graph, the (open) 2--out--neighbourhood of a vertex $a \in [n]$ consists of the vertices reachable from $a$ by a (consecutively oriented) path of at most two edges.  We will say that two directed Hamilton paths in $K_n$ are two--separated if for at least one vertex of $K_n$ the union of its  2--out--neighbourhoods in the two paths contains exactly 4 vertices. Formulating this condition in terms of permutations means that for some $a\in [n]$ its two immediate successors in the two linear orders defining the two respective 
permutations are $axy$ and $avw$ where $x,y,v,w$ are four different elements of $[n].$ This relation is seemingly very similar to the one underlying Theorem \ref{thm:sina}. As a matter of fact, in the permutation language, to two crossing Hamilton paths there correspond two permutations such that for some $a \in [n]$ 
its neighbours are completely different, meaning that we have $xay$ and $vaw$ in the two permutations and $x,y,v,w$ are all different. If all pairs of paths in a Hamilton path family are 
crossing then for any one of the paths only one of the corresponding two permutations can be in the family, and it can be either one of the two that correspond to the same path. Somewhat surprisingly, the cardinality of the  largest family of pairwise two--separated Hamilton paths is of a different order of magnitude from the one in 
Theorem \ref{thm:sina}. In fact, we have the following almost tight result.

\begin{thm}\label{thm:bin}

Let $R(n)$ stand for the maximum cardinality of a family of pairwise two--separated Hamilton paths in $K_n.$ Then

$$\frac{n!}{n^{13}5^{n-10}}\leq R(n) \leq \frac{n!}{2^{\lfloor n/2 \rfloor}}.$$
\end{thm}

\proof

In this proof we will represent our Hamilton paths as permutations. We will say that two permutations are two--separated if they correspond to two--separated Hamilton paths.
We associate with any permutation a sequence of unordered couples of elements of $[n].$ The permutation  $\pi=\pi_1\ldots \pi_n$ is associated with the sequence of unordered and disjoint couples 
$\{\pi_1,\pi_2\}\dots \{\pi_{2i-1}, \pi_{2i}\}\dots \{\pi_{2\lfloor n/2 \rfloor-1}, \pi_{2\lfloor n/2 \rfloor}\}.$ We will call this sequence of couples the couple order of $\pi.$ 
We then partition the set of all permutations of $[n]$ according to their couple order. In each class the permutations it contains have the same couple order and all the permutations with the same couple order belong to this class. We claim that no two permutations in the same class are two--separated. Considering that the number of classes is 
$\frac{n!}{2^{\lfloor n/2 \rfloor}}$, we will thus obtain our upper bound.

In fact, let the permutations $\pi'$ and $\pi''$ have the same couple order and let $xy$ and $vw$ be the two immediate successors of  $a$ in $\pi'$ and $\pi''$, respectively. We will distinguish two cases. In the first case $a$ and $x$ appear in the same couple $\{a, x\}$  of the couple order, while in the second case $a$ is in a different couple from $x$ and $y$, which then 
have to appear together in the subsequent couple. In the first case either $a$ and $x$ appear in the same order $(a,x)$ in both permutations and we are done, or their order is different, but then $y$ must appear in one of the two positions subsequent to $a$ in $\pi'',$ contradicting our hypothesis that  $\pi'$ and $\pi''$ are two--separated. It remains to see what happens if in the couple order of $\pi'$ the number $a$ is in the couple preceding $\{x, y\}.$ This means that in the couple order $a$ is in the same couple with some 
$z \in [n]$ 
different from both $x$ and $y.$ Since our two permutations have the same couple order, we conclude that at least one of $x$ and $y$ is in one of  the two positions immediately 
following $a$ in $\pi''.$

Let us now turn to lower bounding $R(n)$. To this end we will use a greedy algorithm to exhibit a large enough family of permutations with the required property. At each step in the algorithm we choose an arbitrary permutation and eliminate from the choice space all those incompatible with the chosen one. This procedure goes on until the choice space becomes empty.\\
To analyze this algorithm, we need an upper bound on the number of the permutations incompatible with a fixed permutation $\iota$ of $[n]$. Without loss of generality we assume that $\iota$ is the identical permutation mapping each number in itself. Let $\pi$ be a permutation incompatible with $\iota$. We have that 
\begin{equation}\label{uno}
\mbox{for each  $j$, $1\leq j\leq n-2$, if  $\pi_j\not\in \{n,n-1\}$ then $\left\{ \pi_{j+1}\pi_{j+2}\right\} \cap \left\{ \pi_j+1,\pi_j+2\right\}\neq \emptyset$}
\end{equation}
otherwise the two immediate successors of $\pi_j$ in the two permutations $\iota$ and $\pi$ are four different elements of $[n],$ contradicting that $\iota$ and $\pi$ are incompatible.\\
If $\pi_{j+1}\in\left\{ \pi_j+1,\pi_j+2\right\}$ or respectively $\pi_{j+2}\in\left\{ \pi_j+1,\pi_j+2\right\}$ then we say that $\pi_{j+1}$ or respectively $\pi_{j+2}$ is a close enough follower of $\pi_{j}$. In this terminology property (\ref{uno}) means that at least one of $\pi_{j+1}$ and  $\pi_{j+2}$ must be a close enough follower of $\pi_j$.\\
In the following we call {\em run of big jumps} in $\pi$ any maximal (not extendable) sequence $\pi_j\pi_{j+1}\ldots \pi_{j+t} $, $1\leq j\leq n-t$ such that 
$\pi_{j+1}-\pi_j \not \in \{-2, -1, 1,2,3\}$ and $\pi_{k+1}-\pi_k \not \in \{-1,1,2\}$, $j<k<j+t$. We will refer to the position $j$ as  the head of the run and to the positions $k$ of the run with $k\geq j+2$ as tail positions of the run. We call {\em  free positions } of $\pi$ the first position of the entire permutation,  the second position in  any run and the two positions immediately following those of either $n$ or $n-1$. All the other positions will be called constrained. (We will soon explain the meaning of this terminology.)
For any constrained position $j$ not belonging to any run we have $|\pi_j-\pi_{j-1}|\leq 3,$  implying that 
$$\pi_j\in \left\{\pi_{j-1}-2,\pi_{j-1}-1,   \pi_{j-1}+1,\pi_{j-1}+2,\pi_{j-1}+3\right\}$$
and the  same is true if $j$ is the head of a run since in this case either $\pi_j-\pi_{j-1}>3$ or $\pi_j-\pi_{j-1}< -2$ would contradict the maximality of the run.\\
Consider now a constrained position $j$ in the tail of a run. By definition, position $j-2$ also belongs to the run. Thus either $\pi_{j-1}-\pi_{j-2}\geq 3$ or 
$\pi_{j-1}-\pi_{j-2}\leq -2.$ 
This implies that $\pi_{j-1}\not\in\left\{\pi_{j-2}+1,\pi_{j-2}+2\right\}$  and from property (\ref{uno}) we obtain   $\pi_{j}\in\left\{\pi_{j-2}+1,\pi_{j-2}+2\right\}$. 
We have just proved that for any constrained position $j$ of $\pi$ the value $\pi_j$ is close to at least one of  the values $\pi_{j-1}$ and  $\pi_{j-2}$. \\

Let us now consider the permutations incompatible with $\iota$ having exactly $r$ runs of big jumps. In these permutations there are  at most $r+5$ free positions and at most ${n-1 \choose {r+2}}$ possibilities for these free positions (remember that, by definition,  the first position is always a free position and that  the two  free positions depending on $n$ are adjacent to it and the same is true for the two free positions depending on $n-1$). Thus the number of these incompatible permutations is upper bounded by 
\begin{equation}\label{due}
{n-1 \choose {r+2}}n^{r+5}5^{n-(r+5)}. 
\end{equation}
We will prove that in any incompatible permutation there can be at most three runs of big jumps. This in turn implies that  any incompatible permutation can have at most eight free positions and from bound (\ref{due}) we have that the number of permutation incompatible with permutation $\iota$  is bounded above by 
$$4{n-1 \choose {5}}n^85^{n-8}.$$
This means that the greedy algorithm will eliminate at most $n^{13}5^{n-10}$ permutations at each step, yielding a family  with the desired property and containing at least
$$\frac{n!}{n^{13}5^{n-10}}$$
permutations, as claimed for the lower bound.\\
To conclude the proof it remains to show that in a permutation $\pi$ incompatible with $\iota$  there are at most three runs of big jumps.
To this aim we now prove that any run of big jumps of $\pi$ that does not contain neither $n$ nor $n-1$ is a suffix of $\pi$.
The proof is by contradiction. Let us assume that in  $\pi$ there is a run that does not contain neither $n$ nor $n-1$ and yet it is not a suffix. Let $t$ be the last position of such a run.

We distinguish two cases. In the first case, we suppose
$$\pi_{t-1}>\pi_t.$$
By assumption, $\pi_t$ is the last position in the run of big jumps. Hence the last inequality can be satisfied only if 
$$\pi_t\leq \pi_{t-1}-2.$$
Thus $\pi_t$ is not close enough to $\pi_{t-1}$ whence, by property (\ref{uno}), we must have $\pi_{t+1}\in \{\pi_{t-1}+1, \pi_{t-1}+2\}$. (Note that, by our assumption,  $\pi-{t}$ is not the last element in $\pi$). 
But then 
$$\pi_{t+1}-\pi_t\geq 3,$$
contradicting the assumption that $\pi_t$ is in the last position of the run of big jumps.

In the second case we suppose 
$$\pi_{t-1} < \pi_t.$$
This implies 
\begin{equation}\label{con}
\pi_t \geq \pi_{t-1}+3.
\end{equation}
Thus $\pi_t$ is not a close enough follower of  $\pi_{t-1}$ so that once again, by property (\ref{uno}), we must have $\pi_{t+1}\in \{\pi_{t-1}+1, \pi_{t-1}+2\}$.
In particular, we see that $\pi_{t+1}<\pi_t.$ However, since $\pi_{t+1}$ is not part of the run of big jumps, the only possibility left is
\begin{equation}\label{cine}
\pi_{t+1}=\pi_t -1.
\end{equation}
Moreover since  $\pi_t$ is not a close enough follower of  $\pi_{t-1}$ this role must be taken by $\pi_{t+1}$. Hence from (\ref{cine}) and from (\ref{con}) it follows that
\begin{equation}\label{qua}
\pi_t=\pi_{t-1}+3 
\end{equation}
Thus, by
 (\ref{cine}) and (\ref{qua}), we have
\begin{equation}\label{tre}
\pi_{t+1}=\pi_{t-1}+2.
\end{equation}
Let us see what we can say at the light of this about $\pi_{t-2}.$ Note that $\pi_t-\pi_{t-1}=3$ which  implies that position $t-1$ is not  the head of the run of big jumps. This in turn implies that   position $t-2$ is in the run and $\pi_{t-1}-\pi_{t-2}\not\in \{-1,1,2\}$. Hence we  must have $\pi_{t-2}\leq \pi_{t-1}-3$ or $\pi_{t-2}\geq \pi_{t-1}+2$. Moreover,  property (\ref{uno}) implies
$$\{\pi_{t-2}+1, \pi_{t-2}+2\} \cap \{\pi_{t-1}, \pi_t\} \not =\emptyset,$$
whence the only choice left is $\pi_{t-2}=\pi_{t-1}+2$ (remember that $\pi_t=\pi_{t-1}+3$), and  this is impossible by equality (\ref{tre}).

\hfill$\Box$

\section{Problem classes}

We have already mentioned the need for a meta--conjecture to describe for which problems within our framework does the optimal solution have a kernel structure. It seems to us 
that this can happen only in case of one--family problems and it will never be the case if the family $\D$ in a two--family problem is not monotone. (A class $\F$ of graphs is monotone if $G \in \F$ and $G\subseteq H$ imply $H \in \F.$ Note, in particular, that the family $\D_3$ of graphs having a vertex of degree 3 is not monotone.) 
Our problems belong to a class of combinatorial problems introduced in \cite{Siso} where they are referred to as intersection problems.
In the eyes of the authors of \cite{Siso} an intersection problem is in terms of  the maximum cardinality of a family of objects such that the pairwise intersections of these objects have some prescribed property. (Even if, in this paper,  we are considering pairwise unions, we might have considered, equivalently, the intersections of  the complementary objects.) We feel that it is more appropriate to reserve the term intersection problem to the case 
$\F=\D$ when in fact we end up with an intersection problem in disguise, in the sense of \cite{Siso}. In our opinion, an intersection problem is a one--family problem. In such problems one asks for the largest cardinality of a family of objects, typically subgraphs of $K_n$ with a given property such that the intersection of any two members of the family is still in the family. One often considers as intersection problems only those about the maximum cardinality of a family of graphs such that the intersection of any two of them is non--empty. 

A more delicate question is to understand when exactly can the optimal solution of a problem of our kind be obtained through a greedy algorithm. Our two main results can hopefully contribute to a better understanding of this issue. At a first glance one is tempted to believe that for an information--theoretic (genuinely two--family) problem the greedy 
algorithm often gives close-to-optimal constructions. Our question is not precise enough and we do not expect a precise answer at this stage.

\section{Acknowledgement} 

We are grateful to Miki Simonovits for stimulating discussions about the true nature of intersection problems.

\end{document}